\documentclass[a4paper,12pt,draft]{amsart}
\usepackage{amsmath}
\usepackage{amsfonts}
\usepackage{amssymb}
\usepackage[OT4]{fontenc}
\usepackage{fullpage}
\newtheorem{theorem}{Theorem}
\newtheorem{lemma}[theorem]{Lemma}

\newtheorem{proposition}[theorem]{Proposition}
\newtheorem{definition}[theorem]{Definition}
\newtheorem{example}[theorem]{Example}
\numberwithin{theorem}{section}
\numberwithin{equation}{section}
\DeclareMathOperator\dc{{\it dd^c}}
\DeclareMathOperator\co{{\it\mathbb C^n}}
\DeclareMathOperator\bo{{\it\mathbb B^n}}
\DeclareMathOperator\ma{max}
\DeclareMathOperator\gra{\it\nabla\overline {u}_{k,j}\rm}
\DeclareMathOperator\uj{\it u_{k,j}\rm}
\DeclareMathOperator\tuj{\it\overline{u}_{k,j}\rm}

\DeclareMathOperator\we{\wedge}
\DeclareMathOperator\ijk{\it I_k^j\rm}
\DeclareMathOperator\hijk{\chi_{\it I_k^j}\rm}
\DeclareMathOperator\Om{{\it\Omega}}
\title{An inequality for mixed Monge-Amp\`ere measures}
\subjclass[2000]{32U15}
\author{S\l awomir Dinew}

\begin{document}

\maketitle

\begin{abstract}
We generalize an inequality for mixed Monge-Amp\`ere measures from \cite{K1}. We also give an example that shows that our assumptions are sharp. The corresponding result in the setting of compact K\"ahler manifold is also discussed.
\end{abstract}
\section{Introduction}
The complex Monge-Amp\`ere operator acts on a smooth plurisubharmonic function $u$\ defined in $\co$\ by
$$MA(u):=4^nn!\det(\frac{\partial^2u}{\partial z_j\partial\overline{z}_k})d\lambda=(\dc u)^n.$$
For this reason the theory of positive definite matrices has found numerous applications in the study of this operator. We refer to \cite{HJ} for the basic ideas in this vein. For example one can use concavity properties of such matrices to obtain pointwise estimates for the Monge-Amp\`ere operator (see \cite{W1}). In particular the following holds:
\begin{theorem}
 Let $u,\ v$\ be bounded and smooth plurisubharmonic functions such that $(\dc u)^n\geq fd\lambda$, $(\dc v)^n\geq gd\lambda $, where $f$\ and $g$\ are smooth nonnegative functions and $d\lambda$\ is the Lebesgue measure. Then
\begin{equation}\label{1}
(\dc u)^k\we(\dc v)^{n-k}\geq f^{\frac k n}g^{\frac{n-k}{n}}d\lambda
\end{equation}
\begin{equation} \label{2}
(\dc(u+v))^n\geq(f^{\frac 1 n}+g^{\frac 1 n})^nd\lambda
\end{equation}
\end{theorem}
 Those inequalities were first used in \cite{BT1} in order to generalize the action of $MA$\ to nonsmooth functions. Very similar estimates were also used by Demailly in problems in complex geometry (see \cite{De}). Last but not least, such inequalities were used in \cite{K1} together with the comparison principle to produce stability estimates for the complex Monge-Amp\`ere operator.

Of course plurisubharmonic functions are, in general, neither smooth nor even bounded. However the Monge-Amp\`ere operator can still be reasonably defined (in the weak sense) for quite general functions (see \cite{BT1}, \cite{BT2}, \cite{Bl1}, \cite{Bl2}, \cite{Ce1}, \cite{Ce2}). This leads to natural questions, namely whether our inequalities hold in this more general situation. Note also that pluripotential theory deals often with measures singular w.r.t. the Lebesgue measure. An example of such a situation arises when one computes the Monge-Amp\`ere of a maximum of two smooth plurisubharmonic functions. Measures associated to maximal plurisubharmonic functions provide another example. Hence one can pose a question in the following way:

\bf Question:\rm\ \ 
\it Let $u,\ v$\ be plurisubharmonic functions (psh for short) such that their Monge-Amp\`ere masses $(\dc u)^n,\ (\dc v)^n$\ are well defined (in the sense of \cite{Bl2}, or equivalently in the sense of \cite{Ce2}). Let also $\mu$\ be a positive measure and $f,\ g\in L^1(d\mu)$. Suppose that $(\dc u)^n\geq fd\mu,\ (\dc v)^n\geq gd\mu$\ in the sense of measures. Is it true that (again in the sense of measures)\rm
\begin{equation}\label{4}
(\dc u)^k\we(\dc v)^{n-k}\geq f^{\frac k n}g^{\frac{n-k}{n}}d\mu
\end{equation}
\begin{equation}\label{5}
(\dc(u+v))^n\geq(f^{\frac 1 n}+g^{\frac 1 n})^nd\mu
\end{equation}

Some results in this direction are known. Inequalities (\ref{4}) and (\ref{5}) were generalized to nonsmooth $u$\ and $v$, but with the the Lebesgue measure instead of $d\mu$ (see \cite{K1}, \cite{Bl3}). In particular the following holds (see \cite{K1}):
\begin{theorem}\label{thm} Assume that $u$\ and $v$\ are plurisubharmonic and bounded functions in a domain in $\co$. Let also $f,\ g\in\ L^1(d\lambda)$\ are nonnegative functions such that 
$$(\dc u)^n\geq fd\lambda,\ (\dc v)^n\geq gd\lambda$$
holds. Then
$$(\dc u)^k\we(\dc v)^{n-k}\geq f^{\frac k n}g^{\frac{n-k}{n}}d\lambda,\ k=1,\cdots,n$$
\end{theorem}
{\bf Remark.} In \cite{K1} this was stated for continuous $u$\ and $v$, but the proof applies without changes to the more general case. The crucial thing in this result is the Lebesgue measure on the right hand side: this enables us to find a nice approximation sequence of smooth functions and use the smooth version to prove the result.\\

The aim of this note is to study these  inequalities for more general functions. In particular we show that they, in general, do not hold. We also give sharp condition for the measure $\mu$\ ensuring the positive answer. It is very interesting that this condition (vanishing on pluripolar sets) is closely related to uniqueness for the Dirichlet problem for the Monge-Amp\`ere operator. Also the same condition appears when one generalizes the Demailly inequality
\begin{equation}
(\dc\max(u,v))^n\geq\chi_{u\geq v}(\dc u)^n+\chi_{u<v}(\dc v)^n
\end{equation}
(see \cite{KH}).
 Let us state our main result:
\begin{theorem} Let $\mu$\ be a positive measure on $\Om$\ that vanishes on all pluripolar sets. Let $u_1,u_2,\cdots,u_n\in PSH(\Om)$\ be plurisubharmonic functions with well defined Monge-Amp\`ere operator (we can work on a fixed domain $\Om$\ as in \cite{Ce1}, or, instead, we can work with germs of functions as in \cite{Bl2}). Let also $f_i,\ i=1,\cdots,n$\ be nonnegative functions integrable with respect to $\mu$. If
$$(\dc u_i)^n\geq f_id\mu,\ \forall i=1,\cdots,n$$
then
$$\dc u_1\we \dc u_2\we\cdots\we\dc u_n\geq (f_1f_2\cdots f_n)^{\frac 1 n}d\mu$$
\end{theorem}
{\bf Remark.}
To unify the two possible approaches we shall work locally (in a small ball). It will be explained what happens if we have a fixed domain $\Om$\ and all objects are globally defined there. Also, for the sake of brevity, we shall work throughout the note with two functions $u$\ and $v$\ instead of the collection of  $n$\ functions. It will be explained how to get this general case.

{\bf Acknowledgement.} I would like to thank my advisor professor S\l awomir Ko\l odziej for very stimulating discussions.
\section{Definitions}
We shall work mostly in $\co$\ but, since we prove purely local results, everything can be generalized to manifold setting. The corresponding result in the setting of compact K\"ahler manifolds will be discussed in the last section.

As usual $d=\partial +\overline\partial$\ and $d^c:=i(\overline\partial-\partial)$, so $\dc=2i\partial\overline\partial$. The complex Monge-Amp\`ere operator is defined to be
$$(\dc u)^n:=\underbrace{\dc u\we\cdots\we\dc u}_{n-times} $$
This definition at first makes sense only for smooth functions (since we cannot multiply distributions), but classical results in pluripotential theory allow to define this also for bounded plurisubharmonic functions.

A domain $\Om\in\co$\ is called {\it hyperconvex} if it admits a negative plurisubharmonic exhaustion function i.e.
\begin{align*}
u\in PSH(\Om),\ u<0,\ \forall c>0\ \lbrace z\in\Om|u(z)<-c\rbrace\Subset\Om
\end{align*}
A smoothly bounded pseudoconvex domain is called {\it strictly pseudoconvex} if it admits a defining function that is strictly plurisubharmonic in $\overline\Om$.

Let $$\mathcal E_0(\Om)=\lbrace u\in PSH(\Om)\cap L^{\infty}(\Om)|\ u<0, lim_{z\rightarrow\zeta}u(z)=0,\ \forall \zeta\in\partial\Om,\ \int_{\Om}(\dc u)^n\leq\infty\rbrace$$ denote the class of bounded exhaustion functions in a hyperconvex domain. Whenever it is clear with which domain we deed we  shall write simply $\mathcal E_0$\ instead of $\mathcal E_0(\Om)$. 

Let $\mathcal E^p$\ be the class of plurisubharmonic functions consisting of those functions $g$, for which there exists a  sequence $g_j \in\mathcal E_0$\ decreasing towards $g$, such that :

\begin{equation*}    \sup_j\int_{\Omega}(-g_j)^p(dd^cg_j)^n<\infty. 
\end{equation*}
\\
Let $u\in\mathcal E^1$. The quantity 
$$\int_{\Om}-u(\dc u)^n$$
is called {\it $\mathcal E^1$-norm} of $u$. This is not a norm in the classical sense, nevertheless turns out to be a useful tool, as the following result (taken from \cite{Ce1}) shows:
\begin{proposition}\label{e1norm}
Let $u,\ v\in\mathcal E^1(\Om),\ u\leq v$. Then
$$\int_{\Om}-v(\dc v)^n\leq\int_{\Om}-u(\dc u)^n.$$
\end{proposition}
So $\mathcal E^1$-norm of a function is controlled by the $\mathcal E^1$-norm of a smaller function.\\

Now we recall more general Cegrell classes:

 Let $\mathcal F(\Om)$\ be the class consisting of those functions $g$, for which there exists a decreasing sequence $g_j \in\mathcal E_0(\Om)$\ such that :
$$\sup_j\int_{\Omega}(dd^cg_j)^n<\infty.$$
Finally let $\mathcal E$\ be the class of functions $g$\ which locally belong to $\mathcal F$, i.e.
 $$\forall w\in\Om\ \exists U_w\subset\Om,\ U_w-{\rm open}\ \exists g_w\in\mathcal F(\Om):g|_{U_w}=g_w|_{U_w}.$$

These classes are usually called Cegrell classes (with boundary value 0) in the literature. Their main applications come from the fact that one can define uniquely their Monge-Amp\`ere measure in such a way, that it is still continuous for decreasing sequences and moreover the definition is consistent with the classical one for smooth functions (which are basic properties for bounded plurisubharmonic functions).

There are very useful inequalities of H\"older type (due to Cegrell \cite{Ce2}) in the class $\mathcal F$. They show how to control integrals of mixed Monge-Amp\`ere type from above, in terms of $\mathcal E^1$\ norms: 
\begin{align*}
&\int_{\Om}-u(\dc v_1)\we(\dc v_2)\we\cdots\we(\dc v_n)\leq
(\int_{\Om}-u(\dc v_1)^n)^{\frac 1n}\cdots(\int_{\Om}-u(\dc v_n)^n)^{\frac 1n};\\
&\int_{\Om}-u(\dc v_1)\we(\dc v_2)\we\cdots\we(\dc v_n)\leq\\
&\leq
(\int_{\Om}-v_1(\dc v_1)^n)^{\frac 1{n+1}}\cdots(\int_{\Om}-v_n(\dc v_n)^n)^{\frac 1{n+1}}(\int_{\Om}-u(\dc u)^n)^{\frac 1{n+1}}\\
&u\in \mathcal E_0(\Om),\ v_1,\cdots v_n\in\mathcal F(\Om).
\end{align*}
For more details concerning these topics we refer to \cite{Ce1},\cite{Ce2}.

While we shall deal mostly with bounded plurisubharmonic functions in this note, we introduce the teminology of Cegrell classes, since it helps to simplify some proofs and to make other more transparent.

Given a positive Radon measure on a bounded domain $\Om\in\co$\ we define its {\it canonical approximation} (see \cite{K4}):
Let $supp\ \mu$\ be contained in a big cube $I$. Consider a subdivision $\mathcal B_k$ of $I$\ into $3^{2kn}$\ congruent semi open cubes $\ijk,\  j=1,\cdots,3^{2kn}$. It is no loss of generality to assume $\mu(\cup_{\ijk\in\mathcal B_k}\partial(\ijk))=0$\ (otherwise we can shift at each stage the boundaries a bit). Now define
\begin{equation}\label{mik}
\mu_k:=\sum_j\frac{\mu(\ijk\cap\Om)}{dV(\ijk\cap\Om)}\hijk dV,
\end{equation}
where $\hijk$\ is the characteristic function of $\ijk$. Of course $\mu_k$\ is weak* convergent to $\mu$\ and every term $\mu_k$\ has a density in $L^{\infty}$\ with respect to the Lebesgue measure.

Below we list some results that we shall need later on:
\begin{theorem}\label{ceg}
Let $\Om$\ be a smoothly bounded strictly pseudoconvex domain in $\co$\ and let $f\in\mathcal C^{\infty}(\partial\Om)$\ be arbitrary. Let also $\mu$\ be a positive measure on $\Om$\ with finite mass and compact support. Suppose $\mu$\ satisfies the following condition for any $p>\frac{n}{n-1}$:

There is a constant $A=A(p)$\ such that
$$\int_{\Om}(-\phi)^pd\mu\leq A(\int_{\Om}(-\phi)^p(\dc \phi)^n)^{\frac p{n+p}}$$
for any $\phi\in\mathcal E_0$. Then:
\begin{enumerate}
\item There exist $u_k\in\mathcal C(\overline{\Om})$ which solve the Dirichlet problem:
\begin{equation*}
\begin{cases} u_k\in PSH(\Om)\cap C(\overline{\Om})\\
(\dc u_k)^n=\mu_k\\
u_k=f\ on\ \partial\Om
\end{cases}
\end{equation*}
where $\mu_k$\ are the canonical approximants of $\mu$.
\item Define $u:=(limsup_{k\rightarrow\infty}u_k)^{*}$. Then there is a subsequence of $\lbrace u_k\rbrace$\ (which after renumbering we still denote by $\lbrace u_k\rbrace$) such that
$u_k\rightarrow u$\ in $L^1(d\lambda)$.
\item We have for this sequence that
\begin{align*}
&\sup_{k}\int_{\Om}|-u_k|(\dc u_k)^n<\infty\\
&\lim_{k\rightarrow\infty}\int_{\Om}|u-u_k|(\dc u_k)^n=0.
\end{align*}
\end{enumerate}
\end{theorem}
The proof of the first part may be found in \cite{K4}. Other results follow from Theorems 5.1 and 7.7 and Lemmas 5.2, 5.3, 7.8 and 7.9 from \cite{Ce1}. We would like to mention that the condition that all the functions have the same boundary values can be weakened. If, for example,  the boudary values of $u_k$\ form a sequence decreasing towards a bounded upper semicontinuous function (which will be the case we shall use later on) $u:=(limsup_{k\rightarrow\infty}u_k)^{*}$\ still makes perfect sense and, applying line by line the proofs from \cite{Ce1} we get the finiteness and the convergence of the integrals in this situation too.

\begin{theorem}\label{ck}
Suppose $u_j\in PSH(\Om)\cap C(\overline{\Om})$ is a sequence that converges to $u\in PSH(\Om)$\ in $L^1(d\lambda,\Om)$. Suppose also all $u_k$\ (and hence also $u$) have the same continuous boundary values, i.e. $\lim_{z\rightarrow\zeta} u_j(z)=f(\zeta)\ \forall\zeta\in\partial\Om$. If moreover $\lim_{k\rightarrow\infty}\int_{\Om}|u-u_k|(\dc u_k)^n=0$\ then 
$u_k$\ converges to $u$\ in capacity.
\end{theorem}
This result is contained in the proof of  Lemma 2.1 in \cite{CK2}. Again  we can carry the argument from \cite{CK2} also if we let boundary values of $u_j$\ to decrease (to be precise, since in the proof there boundary values are used only to ensure  relative compactness of sets $\lbrace u_j<u-a\rbrace,\ a>0$, it is even better when boundary values of $u_j$\ are bigger than those of $u$).

Below we recall a definition of the envelope of subsolutions:
\begin{definition} Let $\mu$\ be a positive Radon measure defined in a hyperconvex domain $\Om$\ and $f\in\mathcal C(\Om)$. Then
$$U(\mu,f)(z):=\sup\lbrace v(z)|\ v\in PSH\cap L^{\infty}(\Om),\ (\dc v)^n\geq\mu,\ \limsup_{z\rightarrow\zeta}v(z)\leq f(\zeta),\ \forall \zeta\in\Om\rbrace$$
\end{definition}
If we additionally assume that $\Om$\ is strictly pseudoconvex and $\mu$\ is sufficiently regular (in particular any for any positive measure its canonical approximants are regular enough for this purpose), then a result from \cite{Ce1} shows that actually this function solves the Dirichlet problem. In particular $(\dc U(\mu,f))^n=\mu$.\\

Let us also recall some notions in the setting of compact K\"ahler manifold. Note that interactions between pluripotential theory, complex dynamics and differential geometry on such manifolds is a very timely topic and subject of intensive research. Hence we would like to show how our results can be used in this setting either.

Let $X$\ be compact $n$-dimensional K\"ahler manifold equipped with fundamental K\"ahler form $\omega$\ (that is $d$-closed strictly positive globally defined form) given in local coordinates by 
$$\omega={\frac i 2}  \sum_{k,j=1}^n g_{k\overline{j}}dz^k \wedge d\overline z^j$$
 We assume that the metric is normalized so that 
$$\int_X\omega^n=1 .$$ 
Recall that 
\begin{equation*}
 PSH(X,\omega):=\lbrace \phi \in L^1(X,\omega):dd^c\phi \geq -\omega,\ \phi \in\mathcal C^{\uparrow}(X) \rbrace
\end{equation*} 
where $\mathcal C^{\uparrow}(X)$\ denotes the space of upper semicontinuous functions. We call the functions that belong to $PSH(X,\omega)$ $\omega$-plurisubharmonic ($\omega$-psh for short). A very similar notion of {\it admissible} functions (aditionally $\mathcal C^{\infty}$\ smoothness and strict inequality in definiton of $PSH(X,\omega)$\ are imposed) is widely used in differential geometry and has a much longer history.

One would like also to define the Monge-Amp\`ere operator 
$$(\omega_u)^n:=(\omega+dd^c u)^n$$
in this setting. Since locally functions that belong to $PSH(X,\omega)$\ are standard plurisubharmonic function minus a (smooth) potential for the form $\omega$\ one can use classical results from pluripotential theory to define this for bounded $u$. This approach was used in \cite{K1}.
 Recently in a series of articles Guedj and Zeriahi (\cite{GZ1}, \cite{GZ2}) developped intrinsic pluripotential theory on compact K\"ahler manifolds. In particular in \cite{GZ2} they defined Cegrell classes in this setting. We recall below the definition of the class $\mathcal E(X,\omega)$.

For every $u\in PSH(X,\omega)$\ $(\omega+\dc max( u,-j))^n$\ is a well defined probability measure regardless $u$\ is bounded or not. By \cite{GZ2} the sequence of measures $\chi_{\lbrace u>-j\rbrace}(\omega+\dc max( u,-j))^n$\ is always increasing and one defines
$$\mathcal E(X,\omega):=\lbrace u\in PSH(X,\omega)\ |\ \lim_{j\rightarrow\infty}\int_{X}\chi_{\lbrace u>-j\rbrace}(\omega+\dc max( u,-j))^n=1\rbrace.$$ 
These functions might be unbounded, but the integral assumption ensures that there is no room for ''bad things'' to happen near $\lbrace u=-\infty\rbrace$. Then one defines 
$$(\omega+\dc u)^n:=\lim_{j\rightarrow\infty}\chi_{\lbrace u>-j\rbrace}(\omega+\dc max( u,-j))^n.$$
In particular Monge-Amp\`ere measures of functions from $\mathcal E(X,\omega)$\ do not charge pluripolar sets. We refer to \cite{GZ2} for a discussion of that notion.

\section{Proof of the main theorem}

Our idea will be quite similar in spirit to that in the proof of Theorem \ref{thm}. We shall find appropriate sequences $u_j,\ v_j$\ for which  Theorem \ref{thm} holds, and prove that they converge in a suitable way to $u$\ and $v$\ respectively ensuring the weak convergence of $(\dc u_j)^k\we(\dc v_j)^{n-k}$\ towards $(\dc u)^k\we(\dc v)^{n-k}$. This is the point where we need the canonical approximants. The delicate point is that (by an example in \cite{CK1}) the weak star convergence $(\dc u_j)^n$\ towards $(\dc u)^n$\ (even if all the fuctions considered have the same boundary values) does not imply strong enough convergence (that is convergence in capacity) of $u_j$\ towards $u$. So we have to use some special sequences of approximants. On the other hand, taking for example convolutions would be enough for the convergence in capacity but the inequality for the approximants is unclear.

First we consider the case of {\it bounded} $u$\ and $v$. In fact for most applications not requiring the theory of Cegrell classes this case is sufficient.

Now note that the  claimed inequality is a local property, hence it suffices to prove it in a (small) ball $\bo$, such that the functions $u,\ v$\ are defined in a neighbourhood of it. Let $m_j,\ n_j$\ be two sequences of smooth functions on $\partial\bo$, decreasing to $u|_{\partial\bo}$\ and $v|_{\partial\bo}$\ respectively. Let $u_j,\ v_j$\ solve
\begin{equation*}
\begin{cases}
u_j\in PSH(\bo)\cap L^{\infty}(\bo)\\
(\dc u_j)^n=(\dc u)^n_j\\
u_j|_{\partial\bo}=m_j
\end{cases}
\end{equation*}
\begin{equation*}
\begin{cases}
v_j\in PSH(\bo)\cap L^{\infty}(\bo)\\
(\dc v_j)^n=(\dc v)^n_j\\
v_j|_{\partial\bo}=n_j
\end{cases}
\end{equation*}
We recall that $(\dc u)^n_j$\ is the canonical approximation of the measure $(\dc u)^n$.
By Theorem \ref{ceg} such solutions exist.

 Before we proceed we would like to point out some subtleties. If $u$\ and $v$\ were continuous, the sequences $m_j,\ n_j$\ would be redundant (since we can work with merely continuous boundary data as well). This point causes some technical problems in the proof. Also we need here to use the canonical approximants for the measures on the right hand side instead of the measures themselves, for the following reason: The Dirichlet problem
\begin{equation*}
\begin{cases}
u_j\in PSH(\bo)\cap L^{\infty}(\bo)\\
(\dc u_j)^n=(\dc u)^n\\
u_j|_{\partial\bo}=m_j
\end{cases}
\end{equation*}
need not have a solution continuous up to the boundary. 
\begin{proposition}\label{abc}
Let $u_j,\ v_j$\ be as above. Define $u:=(\limsup_{j\rightarrow\infty}u_j)^{*}$,\newline $v:=(\limsup_{j\rightarrow\infty}v_j)^{*}$. Assume also that $u_j$ (resp. $v_j$) tend to $u$\ (resp. $v$) in $L^1(d\lambda)$. Then we have
$$(\dc u_j)^k\we(\dc v_j)^{n-k}\rightharpoonup(\dc u)^k\we(\dc v)^{n-k},\ \forall k\in\lbrace1,\cdots,n\rbrace.$$
\end{proposition}
\begin{proof}
 Note that the inequality 
$$\int_{\Om}(-\phi)^pd\mu\leq A(p)(\int_{\Om}(-\phi)^p(\dc \phi)^n)^{\frac p{n+p}},\ \phi\in\mathcal E_0(\Om)$$
holds for any $p$\ with a constant $A(p)$\ dependent on $p$, if $\mu$\ is the Monge-Amp\`ere measure of a bounded plurisubharmonic function (this follows easily from the theory in \cite{Ce1}).
So, by Theorem \ref{ceg} and the discussion after it, we have (after passing to an appropriate subsequences, which for the sake of brewity, will also be denoted by $u_j,\ v_j$), that
\begin{align*}
&\lim_{k\rightarrow\infty}\int_{\Om}|u-u_k|(\dc u_k)^n=0,\\
&\lim_{k\rightarrow\infty}\int_{\Om}|v-v_k|(\dc v_k)^n=0.
\end{align*}
Now Theorem \ref{ck} (see also the remark after it) give us that
$u_k$, and $v_k$\ converge to $u$\ and $v$\ in capacity.

Now we are almost ready to approximate $(\dc u)^k\we(\dc v)^{n-k}$\ \newline by $(\dc u_j)^k\we(\dc v_j)^{n-k}$. Indeed a theorem of Xing (see \cite{X1}) says that given a locally uniformly bounded sequences of $PSH$\ functions $\lbrace w_k^j\rbrace_{k=1}^{\infty},\ j=1,\cdots,n$\ converging in capacity to $w^j\in PSH$, the corresponding Monge-Amp\`ere measures \newline
$\dc w_k^1\we\cdots\we\dc w_k^n$\ converge weakly to $\dc w^1\we\cdots\we\dc w^n$.

 The only thing that we (apriori) don't know, is whether the sequences $u_j,\ v_j$\ are locally uniformly bounded. However this difficulty can be bypassed by noticing that $u_j,\ v_j$\ are uniformly bounded in $\mathcal E^1$\ norm (see \cite{Ce1} or \cite{K2}): to show this take any
$U\Subset\Om$\ with $cap(U,\bo)<\epsilon$. Then
\begin{align*}
&\int_{U}(\dc u_j)^k\we(\dc v_j)^{n-k}\leq\\
&\leq\int_{\bo}-h_{U,\Om}(\dc (u_j+U(0,-m_j)))^k\we(\dc(v_j+U(0,-n_j)))^{n-k}\leq\\
&\leq (\int_{\bo}-(u_j+U(0,-m_j))(\dc(u_j+U(0,-m_j)))^n)^{\frac k{n+1}}\times\\
&\times(\int_{\bo}-(v_j+U(0,-n_j))(\dc (v_j+U(0,-n_j)))^n)^{\frac{n-k}{n+1}}
(\int_{\bo}-h_{U,\Om}(\dc h_{U,\Om})^n)^{\frac1{n+1}}
\leq\\
&\leq C^{\frac n{n+1}}cap(U,\Om)^{\frac1{n+1}}\leq C^{\frac n{n+1}}\epsilon^{\frac 1{n+1}}
\end{align*}
Where $C$\ is the uniform $\mathcal E^1$\ bound for $u_j,\ v_j$, $h_{U,\Om}$\ is the relative extremal function of $U$ and we have used the H\"older type inequalities, which is legal since $u_j+U(0,-m_j), v_j+U(0,-v_j)$\ belong to $\mathcal E_0$\ (see \cite{Ce2}). The rigorous justification of the uniform $\mathcal E^1$\ bound for the sequences is a bit technical and will be given in Lemma \ref{lem} below.

 Now, again due to uniform $\mathcal E^1$\ bounds, we have
$$cap(\lbrace u_j<-s\rbrace,\Om)\leq\int_{\Om}-\frac{|u_j|}{s}(\dc h_{\lbrace u_j<-s\rbrace,\Om})^n\leq\frac{C}{s}$$
with $C$\ independent of $j$\ and $s$\ (in fact much better estimates can be provided but these are satisfactory for our needs).

Fix $s$\ big enough such that $cap(\lbrace u_j<-s\rbrace,\Om)\leq\epsilon, \forall j$. Then for any test function $\chi$\ we have
\begin{align*}
&|\int_{\bo}\chi((\dc u_j)^k\we(\dc v_j)^{n-k}-(\dc u)^k\we(\dc v)^{n-k})|\leq\\
&\leq|\int_{\lbrace u_j\leq-s\rbrace\cup\lbrace v_j\leq-s\rbrace}\chi((\dc u_j)^k\we(\dc v_j)^{n-k}-(\dc u)^k\we(\dc v)^{n-k})|+\\
&+|\int_{\bo}\chi((\dc max(u_j,-s))^k\we(\dc max(v_j,-s))^{n-k}-(\dc u)^k\we(\dc v)^{n-k})|
\end{align*}
But the first term is arbitrary small by the argument above and the second term tends to $0$\ due to Xing's theorem. So, we obtained the desired result.
\end{proof}
\begin{lemma}\label{lem}
There is an absolute constant $C$\ independent of $j$\ such that
$$\int_{\bo}-(v_j+U(0,-n_j))(\dc (v_j+U(0,-n_j)))^n<C$$
\end{lemma}
\begin{proof}
Consider the function $g_j:=U((\dc v_j)^n,0)$. From \cite{Ce1} we know that  $g_j\in\mathcal E_0$\ and $(\dc g_j)^n=(\dc v_j)^n$. Hence by comparison principle applied to the pair $v_j, g_j+U(0,n_j)$\ we get
$$v_j+U(0,-n_j)\geq g_j+U(0,n_j)+U(0,-n_j).$$
Let $h_j:=U(0,n_j)+U(0,-n_j)$
By inequality from Proposition \ref{e1norm} we get
$$\int_{\bo}-(v_j+U(0,-n_j))(\dc (v_j+U(0,-n_j)))^n\leq\int_{\bo}-(g_j+h_j)(\dc ( g_j+h_j))^n.$$
The last term can be decomposed into a sum of terms of the type
$$\binom{n}{m}\int_{\bo}-(g_j+h_j)(\dc g_j)^m\we(\dc h_j)^{n-m},\ m\in\lbrace 0,\cdots,n\rbrace.$$
 Again by Cegrell inequalities such terms are controlled from above by some product of $\int_{\bo}-g_j(\dc g_j)^n$\ and $\int_{\bo}-h_j(\dc h_j)^n$. But $h_j$\ are uniformly bounded, while 
$$\int_{\bo}-g_j(\dc g_j)^n=\int_{\bo}-g_j(\dc v_j)^n\leq\int_{\bo}-(v_j+U(0,-n_j))(\dc v_j)^n,$$
 by the comparison principle. Now $U(0,-n_j)$\ is uniformly bounded, $(\dc v_j)^n$\ have uniformly bounded total masses, and $\sup_j\int_{\bo}-v_j(\dc v_j)^n$\ is finite by Theorem \ref{ceg}. Hence we have obtained the claimed uniform bound.
\end{proof}
Now we prove our main inequality in this case:
\begin{theorem}
$u,\ v\in PSH(\Om)\cap L^{\infty}(\Om)$ satisfy
$$(\dc u)^n\geq fd\mu,\ (\dc v)^n\geq gd\mu$$
where $\mu$\ is an arbitrary positive measure and $f,\ g\in L^1(\Om,d\mu)$, then
$$(\dc u)^k\we(\dc v)^{n-k}\geq f^{\frac kn}g^{\frac{n-k}n}d\mu$$
\end{theorem}
\begin{proof}
Consider the canonical approximation as in Theorem \ref{abc}. We have that
\begin{align*}
&(\dc u)^k\we(\dc v)^{n-k}=\lim_{j\rightarrow\infty}(\dc u_j)^k\we(\dc v_j)^{n-k}\geq\\
&\limsup_{j\rightarrow\infty}\sum_j\hijk\frac{(\int_{\ijk}(\dc u)^n)^{\frac kn}(\int_{\ijk}(\dc v)^n)^{\frac{n-k}n}}{dV(\ijk)}dV\geq\\
&\geq
\limsup_{j\rightarrow\infty}\sum_j\hijk\frac
{(\int_{\ijk}fd\mu)^{\frac kn}(\int_{\ijk}gd\mu)^{\frac{n-k}n}}{dV(\ijk)}dV\geq\\
&\geq\limsup_{j\rightarrow\infty}\sum_j\hijk\frac
{(\int_{\ijk}f^{\frac kn}g^{\frac{n-k}n}d\mu)}{dV(\ijk)}dV=f^{\frac kn}g^{\frac{n-k}n}d\mu
\end{align*}
where we have used Theorem \ref{thm} and the H\"older inequality.
\end{proof}
The case of $n$\ different functions instead of just two goes in the same way. The only difference is that we must use the generalised H\"older inequality (for $n$\ functions) instead of the classical one that we used above.

This result can be generalised to unbounded plurisubharmonic functions. We show below that our inequality remains true provided $\mu$\ does not charge pluripolar sets. As the example in the next section shows the result is sharp. Since this is a purely local result we state it in terms of the Cegrell classes in a hyperconvex domain.
\begin{theorem}
$u,\ v\in PSH(\Om)\cap \mathcal E(\Om)$ satisfy
$$(\dc u)^n\geq fd\mu,\ (\dc v)^n\geq gd\mu$$
assume moreover that $\mu$\ does not charge pluripolar sets. Then we have the same conclusion as in the above theorem.
\end{theorem}
\begin{proof}
We recall the following known inequality which is a special case of Demailly's inequality (see for example \cite{KH}):
\begin{equation}\label{56}
(\dc \max(u,-j))^n\geq\chi_{\lbrace u>-j\rbrace}(\dc u)^n
\end{equation}
for every $u$\ in $\mathcal E(\Om)$.
By the monotone convergence and the result in the bounded case we obtain
\begin{align*}
&(\dc u)^k\we(\dc v)^{n-k}=lim_{j\rightarrow\infty}(\dc\max( u,-j))^k\we(\dc\max(v,-j))^{n-k}\geq\\ 
&\geq limsup_{j\rightarrow\infty}
(\chi_{\lbrace u>-j\rbrace}f)^{\frac kn}(\chi_{\lbrace v>-j\rbrace}g)^{\frac {n-k}n}\mu
\end{align*}
the last term converges to $f^{\frac kn}g^{\frac {n-k}n}\mu$\ (because $\mu$\ does not charge the pluripolar set \newline $\lbrace u=-\infty\rbrace\cup\lbrace v=-\infty\rbrace$), which proves the claim.
\end{proof}

\section{A counterexample}
All the results above suggest that the inequalities might be true for every measure $\mu$. Below we shall present a counterexample which shows that this is not the case. Of course such a measure must charge some pluripolar set.

 This example is borrowed from Wiklund's paper \cite{W2}, where these functions were used in a different context:
\begin{example}
 Let\newline $u_k=\ma\lbrace\frac{1}{k}\log|z_1|,k^2\log|z_2|\rbrace,\ \ v_k=\ma\lbrace\frac{1}{k}\log|z_2|,k^2\log|z_1|\rbrace$. Then
$$(\dc u_k)^2=(2\pi)^2\frac k 2\delta_0,\ 
(\dc v_k)^2=(2\pi)^2\frac k 2\delta_0$$
but
$$\dc u_k\we\dc v_k=(2\pi)^2\frac1{2k^2}\delta_0$$
$$(\dc(u_k+v_k))^2=(2\pi)^2(k+\frac{1}{k})\delta_0$$
where $\delta_0$\ is the Dirac delta. In particular inequalities (\ref{4}) and (\ref{5}) both fail in this case. One can also show that Demailly inequality fails for these functions, so the result in \cite{KH} is sharp too.
\end{example}
\begin{proof}
First let us compute $(\dc u_k)^2$. This is already known (see \cite{W2}) but first, we would like to outline the proof and secondly, we shall need that method later. Note that our function is pluriharmonic except on the (real) hypersurface $|z_1|=|z_2|^{k^3}$. It is enough to compute $(\dc\uj)^2$\ with $\uj:=\ma\lbrace u_k,-j\rbrace$, since by a result of B\l ocki such a sequence of measures weakly converege to $(\dc u_k)^2$\ for this particular choice of $u_k$\ (because this function belongs to $W_{loc}^{1,2}$, or equivalently it belongs to the Cegrell class $\mathcal E$, see \cite{Bl1}, \cite{Ce1}). Now proceeding as in \cite{Bl1} we use the change of the variable
$$(x,y)\longrightarrow(\log|z_1|,\log|z_2|)$$
to confirm that
$$\int_{\mathbb D^2}(\dc\uj)^2=\int_{x\leq0, y\leq 0}MA(\tuj)$$
where $MA$\ is the real Monge-Amp\`ere operator and $\tuj(x,y):=\ma\lbrace\frac 1 kx,\ k^2y,\ -j\rbrace$. By Alexandrov's theorem (see \cite{Al}) the latter integral is equal to the volume of the gradient image, i.e.
$$\int_{x\leq0, y\leq 0}MA(\tuj)=\lambda(\gra(\lbrace x\leq0, y\leq 0\rbrace)),\ \gra(E):=\cup_{w\in E}\gra(w)$$
where $\gra(w):=\lbrace t\in\mathbb R^n|\uj(w)+<s-w,t>\leq\uj(s),\ \forall s\in Dom\uj\rbrace$.

At points where $\uj$\ is smooth $\gra(w)$\ is a singleton set consisting of the usual gradient of $\uj$, while at non-smooth points usually $\gra$\ is not a singleton. Hence at points where $\uj$\ is smooth and equal to $\frac 1 k x$\ we get that $\gra(w)=\lbrace(\frac 1 k,0)\rbrace$. Analogously in the two other  smooth regions $\uj=k^2y$\ and $\uj=-j$\ we get that $\gra(w)$\ is equal to $\lbrace(0,k^2)\rbrace$\ and $\lbrace(0,0)\rbrace$, respectively. Note that the Lebesgue measure of the gradient image for this set is $0$. Let now $w$\ is a point where (for example) $\frac1 kx=k^2y>-j$. Then one easily computes the gradient image to be the line segment joining $(\frac 1 k,0)$\ and $(0,k^2)$. Analogously for the other points where two of the three functions considered in the maximum coincide the gradient image is a line segment joining the corresponding endpoints. Finally at the point$(-kj,\frac{-j}{k^3})$\ (all three functions coincide), the gradient image will be the full triangle with vertices $(\frac 1 k,0),\ (0,k^2),\ (0,0)$.

The analysis above shows us that the total mass of $(\dc\uj)^2$\ over the unit bidisc is equal to  $(2\pi)^2\frac k 2$. Also if we fix a set $U\subset\mathbb D^2$\ disjoint from the origin, its logarithmic image would not contain $(-kj,\frac{-j}{k^3})$\ for $j$\ large, hence the gradient image of that set will have zero Lebesgue measure. This shows that $(\dc u_k)^2$\ is concentrated at the origin and since we know the total mass we find that $(\dc u_k)^2=(2\pi)^2\frac k 2\delta_0$. Note that $(\dc v_k)^2=(2\pi)^2\frac k 2\delta_0$\ by symmetry.

What is left is to compute $\dc u_k\we\dc v_k$. But note that
$$2dd^c u_k\we\dc v_k=(\dc(u_k+v_k))^2-(\dc u_k)^2-(\dc v_k)^2$$
so all we need to do is to compute $(\dc(u_k+v_k))^2$. Note that
$$u_k+v_k=\ma\lbrace(k^2+\frac 1 k)\log|z_1|,\ \frac 1 k \log|z_1z_2|,\ (k^2+\frac 1 k)\log|z_1|\rbrace.$$
 Arguing in the same way the gradient image of $\ma\lbrace\overline{u_k+v_k},\ -j\rbrace$ is the (obtuse) rectangle with vertices (clockwise) $(0,k^2+\frac1 k),\ (\frac 1k,\frac1 k),\ (k^2+\frac1 k,0),\ (0,0) $ which has volume $k+\frac 1{k^2}$. Hence as above $(\dc(u_k+v_k))^2=(2\pi)^2(k+\frac1{k^2})\delta_0$, and finally  $\dc u_k\we\dc v_k=(2\pi)^2\frac1{2k^2}\delta_0$.
\end{proof}
\section{the K\"ahler manifold case}
In this section we prove analogous inequalities on a compact K\"ahler manifold.
\begin{theorem}
Let $u,\ v\in\mathcal E(X,\omega)$\ be $\omega$-plurisubharmonic functions on $X$, and let $\mu$\ be a positive measure that does not charge pluripolar sets and $f,\ g\in L^1(d\mu)$. If
$$(\omega+\dc u)^n\geq fd\mu,\ (\omega+\dc v)^n\geq gd\mu$$
as measures, then 
$$(\omega+\dc u)^k\we(\omega +\dc v)^{n-k}\geq f^{\frac kn}g^{\frac{n-k}n}d\mu,\ \forall k\in\lbrace 1,\cdots,n-1\rbrace$$
\end{theorem}
\begin{proof} By the definition of $\mathcal E(X,\omega)$\ it is enough to work with bounded approximants (on the set $\lbrace u>-j\rbrace\cap\lbrace v>-j\rbrace$). Also since again the problem is local it is enough to prove such an inequality in a ball contained in a coordinate chart. But then it is equivalent to the result we have already obtained.
\end{proof}
Finally we would like to mention a slight generalisation of the result above. Namely it is redundant to impose that both functions belong to the same space $\mathcal E(X,\omega)$. We have the following proposition:
\begin{proposition}
Let $u\in\mathcal E(X,\alpha), v\in\mathcal E(X,\beta)$\ be $\omega$-plurisubharmonic functions on $X$\ ($\alpha$\ and $\beta$\ are possibly different K\"ahler forms on $X$), $\mu$\ be a positive measure that does not charge pluripolar sets and $f,\ g\in L^1(d\mu)$. If
$$(\alpha+\dc u)^n\geq fd\mu,\ (\beta+\dc v)^n\geq gd\mu$$
as measures, then 
$$(\alpha+\dc u)^k\we(\beta +\dc v)^{n-k}\geq f^{\frac kn}g^{\frac{n-k}n}d\mu,\ \forall k\in\lbrace 1,\cdots,n-1\rbrace$$
\end{proposition}
\begin{proof} The proof above applies, since the argument for the transition from a ball in a chart to a ball in $\co$\ is independent of the potential for the K\"ahler form.
\end{proof}

S\l awomir Dinew\\
Institute of Mathematics\\
Jagiellonian University\\
ul. Reymonta 4\\
30-059 Krak\'ow\\
Poland\\
\tt slawomir.dinew@im.uj.edu.pl
\end{document}